\documentclass[11pt]{amsart}
\usepackage[utf8]{inputenc}
\usepackage{amsmath} 
\usepackage{amssymb} 
\usepackage{amsthm,amsrefs} 
\usepackage{amsfonts} 
\usepackage{mathtools}
\usepackage{pgfplots}

\usepackage{cases}

\usepackage{enumerate}

\usetikzlibrary{positioning}
\usetikzlibrary{angles,quotes}

\pgfplotsset{compat=1.18, width=12cm}
\usepackage[hidelinks,colorlinks=true,linkcolor=blue,citecolor=blue]{hyperref}
\usepackage[margin= 1in]{geometry}
\usepackage{tabularx}

\usepackage{graphicx}

\renewcommand{\phi}{\varphi}

\usepackage[capitalise]{cleveref}

\usepackage{pythonhighlight}

\newcommand{\ol}[1]{\overline{#1}}
\newcommand{\mb}[1]{\mathbb{#1}}
\newcommand{\mc}[1]{\mathcal{#1}}

\newcommand*{\rom}[1]{\expandafter\romannumeral #1}

\newcommand{\mm}{\mathbb{R}}
\newcommand{\m}{t}
\newcommand{\mmone}{t-1}

\newcommand{\cK}{{\mathcal{K}}}

\newtheorem{theorem}{Theorem}
\newtheorem{corollary}[theorem]{Corollary}

\newtheorem{lemma}[theorem]{Lemma}

\newtheorem{proposition}[theorem]{Proposition}
\newtheorem{remark}[theorem]{Remark}

\title{Hadwiger's conjecture for cap bodies}

\author{Andrii Arman}
\address{Department of Mathematics, University of Manitoba, Winnipeg, MB, R3T 2N2, Canada}
\email{andrew0arman@gmail.com}
\thanks{}

\author{Jaskaran Singh Kaire}
\email{singhj82@myumanitoba.ca}
\thanks{The second author was supported in part by NSERC of Canada Discovery Grant RGPIN-2020-05357.}

\author{Andriy Prymak}
\email{prymak@gmail.com}
\thanks{The third author was supported by NSERC of Canada Discovery Grant RGPIN-2020-05357.}



\begin{document}
\begin{abstract}
Hadwiger's covering conjecture states that every $n$-dimensional convex body can be covered by at most $2^n$ of its smaller positive homothetic translates, with $2^n$ copies required only for affine images of the $n$-cube. 

In this note, we confirm Hadwiger's conjecture for the class of cap bodies in all dimensions, bridging recently established cases of $n=3$ and large $n$. The proof uses probabilistic
techniques and, for $4\le n \le 15$, computer-assisted linear programming.

\end{abstract}

\maketitle

\section{Introduction}

Hadwiger~\cite{H} asked the following question in 1957: For $n\geq 3$, what is the smallest number $H(n)$ such that every $n$-dimensional convex body can be covered by the union of at most $H(n)$ translates of the body's interior? An equivalent formulation in terms of illumination was offered by Boltyanski~\cite{B} in 1960. For a convex body $K$, a direction (unit vector) $v$ illuminates a point $x$ on the boundary $\partial K$ of $K$ if the ray $\{x+vt \;:\; t\geq 0\}$ has nonempty intersection with the interior of $K$. The set of directions  $\{v_i\}_{i=1}^k$ is said to illuminate $K$ if every point of $\partial K$ is illuminated by some $v_i$. The illumination number $I(K)$ of $K$ is the smallest $k$ for which $K$ can be illuminated by $k$ directions. For a collection $\mathcal{K}$ of convex bodies, we define the illumination number as $I(\mathcal{K}):=\max\{I(K):K\in\mathcal{K}\}$. Hadwiger's covering number $H(n)$ is then the illumination number of the class of all $n$-dimensional convex bodies.

The Hadwiger covering conjecture, also known as Levi--Hadwiger(--Gohberg--Markus) conjecture, or Hadwiger--Boltyanski conjecture, states that $H(n)=2^n$ with $I(K)=2^n$ if and only if $K$ is an affine copy of an $n$-cube.
The conjecture attracted many mathematicians, but remains unsolved for all $n\ge 3$. For further details and various partial cases, an interested reader is referred to the survey~\cite{BK}, and more recent works~\cite{Su-Vr}, \cite{ABP-hadwiger}.

The goal of this paper is to complete the confirmation of Hadwiger's conjecture for the class of cap bodies in all dimensions. Following
the terminology of~\cite{BIS}, the convex hull of a ball and an external point is called a spike (or cap). The union of finitely many spikes of a ball is called a cap body if it is a convex set. Let $\mathcal{K}^n_c$ denote the class of all $n$-dimensional cap bodies, while $\mathcal{K}^n_{c,s}$ and $\mathcal{K}^n_{c,us}$ denote the subclasses of $\mathcal{K}^n_{c}$ of all centrally symmetric and unconditionally symmetric (about every coordinate hyperplane) cap bodies, respectively.

Cap bodies were considered by Nasz\'odi in~\cite{Na} (called spiky balls there) as examples of convex bodies which can be ``very close'' to the ball while having exponentially large (in the dimension $n$) illumination numbers. Ivanov and Strachan~\cite{IS} showed that $I(\cK^3_{c,s})=6$ and $I(\cK^4_{c,us})=8$. Bezdek, Ivanov and Strachan~\cite{BIS} proved that $I(\mathcal{K}^n_{c,s})<2^n$ for $n=3,4,9$ and all $n\ge19$, thereby confirming the conjecture for centrally symmetric cap bodies in these dimensions. They used a reduction to symmetric coverings of the sphere with spherical caps of radius $\pi/4$. For large $n$, the estimate obtained in~\cite{BIS} is actually $I(\mathcal{K}^n_{c,s})<(\sqrt{2}+o(1))^n$. This was improved in~\cite{ABPR} using a different method combining spherical coverings and packings, which also allowed us to get rid of the symmetry assumption. Namely, \cite{ABPR}*{Th.~5} implies $I(\mathcal{K}^n_c)<1.19851^n$ for all $n\ge n_0$, where $n_0$ is implicit, confirming the conjecture for all (not necessarily centrally symmetric) cap bodies in high dimensions. In the three-dimensional case, we recently proved~\cite{AKP-3d} that $I(\cK^3_c)=6$. The proof is based on choosing 4 of the illuminating directions as vertices of a randomly rotated regular simplex, and showing, by reducing the problem to an integer linear program and using computer assistance, that for any given cap body $K$, with positive probability at most 2 caps remain unilluminated, so an additional 2 directions bring the illumination number of $K$ to at most 6.

In this paper, we extend and generalize the techniques of~\cite{AKP-3d} to show that $I(\cK^n_c)<2^n$ for all $n$, and thus, we confirm Hadwiger's conjecture for cap bodies in all dimensions. More precisely, for small/moderate dimensions, we illuminate by the directions of several randomly rotated regular simplices and/or several randomly rotated cross-polytopes, and estimate, by linear programming and computer assistance, the expected number of unilluminated caps. We obtain:
\begin{theorem}
    \label{thm:small}
    $I(\cK^n_c)<2^n$ for all $4\le n\le 15$, moreover, $I(\cK^4_c)\le 11$, $I(\cK^5_c)\le 17$, $I(\cK^6_c)\le 27$, $I(\cK^7_c)\le 44$, $I(\cK^8_c)\le 66$. (See \cref{tbl} in \cref{sec:small} for more details.) 
\end{theorem}
One of the key ingredients for the estimates of the expected
number of unilluminated caps is a bound on the intersection of two equal hyperspherical caps, given in \cref{lem:area-inter-equal}. 

For large dimensions, we use an appropriate number of random rotations of the cross-polytope and direct estimates (no computer assistance required) to prove the following:
\begin{theorem}
    \label{thm:explicit}
    For all $n\ge 9$
    \begin{equation}
        \label{eqn:explicit}
        I(\cK^n_{c})<n+7+(\sqrt{2})^n\sqrt{\pi n}(1+\tfrac3n+\tfrac12n\ln\tfrac n2 +\ln\tfrac2{\sqrt{n-1}}).
    \end{equation}
    In particular, $I(\cK^n_{c})<2^n$ for all $n\ge 13$.
\end{theorem}
While asymptotically~\eqref{eqn:explicit} gives $I(\mathcal{K}^n_{c})<(\sqrt{2}+o(1))^n$, which is weaker than in \cite{ABPR}*{Th.~5}, our bound here is explicit and provides the desired bound for all $n\ge 13$. 

Now \cref{thm:small,thm:explicit}, together with $I(\cK^3_c)=6$ (\cite{AKP-3d}*{Th.~2}) imply that the Hadwiger conjecture is valid for all cap bodies:
\begin{corollary}
    $I(\cK^n_c)<2^n$ for $n\ge 3$.
\end{corollary}
Note that $I(\cK^2_c)=4$, as the upper bound is valid for all convex bodies (see ~\cite{Le}), while $\cK^2_c$ contains a square, which has illumination number $4$. At the same time, $\cK^n_c$ for $n\ge 3$ does not contain any affine copy of $n$-cube.

We gather all preliminaries in \cref{sec:prel}. \cref{thm:small,thm:explicit} are proved in \cref{sec:small,sec:large}, respectively.


\section{Preliminaries}\label{sec:prel}

\subsection{Cap bodies and illumination.}
Let $\mm{}^n$ denote the $n$-dimensional Euclidean space with the Euclidean inner product $\langle \cdot,\cdot \rangle$ and the norm $\|\cdot\|$. The unit sphere and the unit ball centred at the origin are denoted by $\mb{S}^{n-1}:=\{x\in \mm{}^n: \|x\|=1\}$ and $\mb{B}^{n}:=\{x\in \mm{}^n: \|x\|\leq 1\}$, respectively. For any points $x,y\in \mb{S}^{n-1}$, the geodesic distance between them is defined by $\theta(x,y):=\arccos\langle x,y\rangle$. For $\xi \in \mb{S}^{n-1}$, define the open and closed spherical caps on $\mb{S}^{n-1}$ centred at $\xi$ of radius $\varphi$ by $C(\xi,\varphi):=\{y \in \mb{S}^{n-1}: \langle \xi, y \rangle> \cos\varphi\}$, $C[\xi,\varphi]:=\{y \in \mb{S}^{n-1} : \langle \xi, y \rangle\geq \cos\varphi\}$.

Let $\text{conv}(X)$ be the convex hull of the set $X$. $K$ is a convex body in $\mm{}^n$ if it is a compact convex set with non-empty interior. A convex body $K$ is called a cap body 
if and only if  $$K=\bigcup_{i=1}^m\text{conv}(\{x_i\} \cup \mb{B}^{n})$$ for some points $x_i\in \mm{}^n\setminus \mb{B}^{n}$, $1\le i\le m$, which are called vertices of $K$. For a given vertex $x_i$, the corresponding base cap (or simply cap, with a slight abuse of
terminology) is defined to be the set 
\[
S_i:={\mathrm{cl}\,}\left({\text{conv}(\{x_i\}\cup \mb{B}^{n})\setminus \mb{B}^n}\right)\cap \mb{S}^{n-1},
\]
where ${\mathrm{cl}\,}{(\cdot)}$ denotes the closure of the set.
Note that $S_i = C[\widehat x_i,\varphi_i]$, where $\widehat{x}_i:=\frac{x_i}{\|x_i\|}$ is the centre of the cap, and $\varphi_i=\arccos{}\frac{1}{\|x_i\|}$ is the radius of the cap. It is an easy observation that the base caps always have acute radius, i.e., $\varphi_i<\pi/2$ for $1\le i\le m$. Furthermore, observe that the convexity of $K$ implies that 
\begin{equation}
    \label{eqn:caps-not-overlap}
C(\widehat x_\alpha,\varphi_\alpha)\cap C(\widehat x_\beta,\varphi_\beta)=\emptyset \quad\text{for any distinct}\quad  x_\alpha, x_\beta\in \{x_i\}_{i=1}^m.
\end{equation}


We use the following proposition from~\cite{BIS,ABPR}:
\begin{proposition}
    \label{prop:ABPR} A cap body $K$ with vertices $\{x_i\}_{i=1}^m$ is illuminated by the directions $\{v_j\}^k_{j=1} \subset \mb{S}^{n-1}$ if:
    \begin{enumerate}[i)]
    \item $C\left(-\widehat{x}_i,\frac{\pi}{2}-\varphi_i\right)\cap \{v_j\}^k_{j=1} \neq \emptyset$ for each $i$, $1\leq i\leq m$,
    \item positive hull of $\{v_i\}^k_{i=1}$ is $\mm{}^n$, i.e. for all $x
    \in \mm{}^n$ there are positive $c_1,\ldots,c_k$ such that $x=c_1v_1+\cdots+c_kv_k$. 
\end{enumerate}
\end{proposition}

\subsection{Area of spherical cap.}
Let $\sigma$ be the probabilistic spherical measure on $\mb{S}^{n-1}$. It is well-known (e.g. see~\cite{Li}) that  $\sigma(C[x,\theta])=\frac{1}{2}I_{\sin^2\theta}\left(\frac{n-1}{2},\frac{1}{2}\right)$
for $0\le\theta\le\pi/2$, where $$I_z(a,b)=\frac{1}{B(a,b)}\int_0^zu^{a-1}(1-u)^{b-1}\,du$$ is the regularized ($I_1(a,b)=1$) incomplete beta function. We will also use the estimate (see, e.g.~\cite{BW}, proof of Cor.~3.2(i), (iii))
\begin{equation}
    \label{eqn:BW}
    \frac{\sin^{n-1}\theta}{\sqrt{2\pi n}} \le \sigma(C[x,\theta])\le \frac{\sin^{n-1}\theta}{\sqrt{2\pi (n-1)} \cdot \cos \theta},
\end{equation}
where the upper bound is valid for $\theta<\pi/2$.

\subsection{Area of intersection of two equal spherical caps.} We need the following technical result.
\begin{lemma}
    \label{lem:area-inter-equal}
    Suppose $n\geq 4$, $0\le\beta\le\alpha\le \pi/2$, and $x,y\in\mb{S}^{n-1}$ with $\theta(x,y)=2\beta$. Then
    \begin{equation}
        \label{eqn:area-inter-equal}
        \sigma(C[x,\alpha]\cap C[y,\alpha])=    \frac{n-2}{\pi}\int_\frac{\cos\alpha}{\cos\beta}^1(1-r^2)^\frac{n-4}{2}\left(\arccos\Bigl(\frac{\cos\alpha}{r}\Bigr)-\beta\right)rdr
        =: A_n(\alpha,\beta).
    \end{equation}
\end{lemma}
\begin{proof}
Suppose that $\mc{D}\subset \{(x_1, x_2) : x_1^2+x_2^2 \le 1\}$ is a measurable subset of the unit disk and $C_{n,\mc{D}}=\{(x_1,x_2,\dots, x_n)\in \mb{S}^{n-1} : (x_1,x_2)\in \mc{D}\}$. First, we will calculate $\sigma(C_{n,\mc{D}})$.

We use standard spherical coordinates
\[
x_i=\cos \phi_i\prod_{1\le j<i}\sin \phi_j , \quad 1\le i<n,\quad x_n=\prod_{j=1}^{n-1}\sin \phi_j, \quad \phi_1,\dots,\phi_{n-2}\in[0,\pi], \quad \phi_{n-1}\in[0,2\pi),
\]
and by ${C_{n,\mc{D}}'}$ and ${\mc{D}'}$ denote ${C_{n,\mc{D}}}$ and ${\mc{D}}$ in $(\phi_1,\dots,\phi_{n-1})$ and $(\phi_1,\phi_2)$ coordinates, respectively. If $|\cdot|$ denotes the Lebesgue surface measure on the unit sphere, then  
\begin{align}
|\mb{S}^{n-1}|\,\sigma(C_{n,\mc{D}}) &= \int_{{C_{n,\mc{D}}'}} \sin^{n-2} \phi_1 \sin^{n-3} \phi_2 \cdots \sin \phi_{n-2}\, d\phi_1\dots d\phi_{n-1} \nonumber \\
&=|\mb{S}^{n-3}|\iint_{{\mc{D}'}} \sin^{n-2}\phi_1 \sin^{n-3}\phi_2 \, d\phi_1d\phi_2 \, \nonumber.
\end{align}

Going back to Cartesian coordinates, $dx_1dx_2=\sin^2\phi_1\sin\phi_2\,d\phi_1 d\phi_2$ and $\sin \phi_1\sin \phi_2=\sqrt{1-x_1^2-x_2^2}$, so using $|\mb{S}^{n-1}|=2\pi^{\frac{n}{2}}/\Gamma(\frac{n}{2})$,
\begin{align}
\,\sigma(C_{n,\mc{D}}) &= \frac{n-2}{2\pi} \iint_{{\mc{D}}} (1-x_1^2-x_2^2)^{\tfrac{n-4}2} \, dx_1 dx_2 \,. \label{eqn:measure C_}
\end{align}

Now we apply~\eqref{eqn:measure C_} to compute  $\sigma(C[x,\alpha]\cap C[y,\alpha])$. Without loss of generality, we may assume that $x=(\cos\beta,\sin\beta,0,\ldots,0)$, $y=(\cos\beta,-\sin\beta,0,\ldots,0)$, in which case  $\sigma(C[x,\alpha]\cap C[y,\alpha])$ is symmetric with respect to $x_2=0$ hyperplane. So $\sigma(C[x,\alpha]\cap C[y,\alpha])=2\sigma(C_{n,\mc{D}_{\alpha,\beta}})$, where \[
\mc{D}_{\alpha,\beta}=\{(x_1,x_2) : x_1^2+x_2^2\le 1, x_2\ge 0, (x_1,x_2)\cdot(\cos\beta,-\sin\beta)\ge \cos\alpha\}.
\]  
In polar coordinates $x_1=r\cos\theta$, $x_2=r\sin\theta$,
we have
\[
\mc{D}^\prime_{\alpha,\beta}=\left\{(r,\theta) : 0\le \theta\le \arccos\left(\frac{\cos\alpha}{r}\right)-\beta, \frac{\cos\alpha}{\cos\beta}\le r \le 1\right\}.
\]
Now~\eqref{eqn:area-inter-equal} readily follows from~\eqref{eqn:measure C_}.
\end{proof}


Clearly, the area of a single cap can be obtained as $\sigma(C[x,\alpha])=A_n(\alpha,0)$, so we will use the notation $A_n(\alpha):=A_n(\alpha,0)$.

\subsection{Union of equal caps centred at the vertices of a simplex or a cross-polytope.} Let $S_n$ be the $n+1$ vertices of a regular simplex inscribed into $\mb{S}^{n-1}$, and let $C_n$ be the vertices of the cross-polytope, i.e. points with coordinates $(\pm1,0,\dots,0)$ and permutations thereof.  For a discrete subset $X\subset\mb{S}^{n-1}$ and $\theta\in(0,\tfrac\pi2]$, denote $C[X,\theta]:=\bigcup_{x\in X}C[x,\theta]$. We need values or lower bounds of $\sigma(C[X,\theta])$ when $X$ is either $S_n$ or $C_n$. We assume $n\ge 3$ below in this subsection.

\begin{lemma}
    \label{lem:Sprob_comp}
Suppose $\psi^{\scriptscriptstyle \triangle}_{n}\in(\tfrac12\arccos(-\tfrac1n), \arccos\tfrac1n]$ is arbitrary, and
\begin{equation*}
    B^{\scriptscriptstyle\triangle}_{n}(\theta):=(n+1)A_n(\theta) - \tfrac{n(n+1)}{2} A_n\!\left(\theta,\tfrac12\arccos(-\tfrac1n)\right), \quad\tfrac12\arccos(-\tfrac1n)< \theta\le\tfrac\pi2. 
\end{equation*}
We have:
\begin{equation*}
\sigma(C[S_n,\theta]) \ 
\begin{cases}
=(n+1)A_n(\theta), & 0<\theta\le \tfrac12\arccos(-\tfrac1n), \\
\ge B^{\scriptscriptstyle \triangle}_{n}(\theta), & \tfrac12\arccos(-\tfrac1n)< \theta \le\psi^{{\scriptscriptstyle \triangle}}_{n},\\
\ge B^{\scriptscriptstyle \triangle}_{n}(\psi^{{\scriptscriptstyle \triangle}}_{n}), & \psi^{{\scriptscriptstyle \triangle}}_{n}< \theta \le\arccos\tfrac1n, \\
=1, & \arccos\tfrac1n \leq \theta\le\tfrac\pi2.
\end{cases}
\end{equation*}
    \end{lemma}
    \begin{proof}
        The geodesic distance between any two distinct points in $S_n$ is $\arccos(-\tfrac1n)$. Therefore,  when $0<\theta\le \tfrac12\arccos(-\tfrac1n)$ the set $C[S_n,\theta]$ consists of $n+1$ disjoint caps of radius $\theta$, and so has measure $(n+1)A_n(\theta)$. When $\tfrac12\arccos(-\tfrac1n)< \theta\le \arccos\tfrac1n$, any two caps in $C[S_n,\theta]$ intersect, so $\sigma(C[S_n,\theta])\ge B^{\scriptscriptstyle \triangle}_{n}(\theta)$ by the Bonferroni's inequality and~\cref{lem:area-inter-equal}. (In fact, $\sigma(C[S_n,\theta])=B^{\scriptscriptstyle \triangle}_{n}(\theta)$ for $\theta\le \arccos{\sqrt{\frac{n-2}{3n}}}$.) 
        
        The case $\psi^{{\scriptscriptstyle \triangle}}_{n}< \theta \le\arccos\tfrac1n$ is merely the consequence of monotonicity of $\sigma(C[S_n,\theta])$ w.r.t. $\theta$ and the already established bound. Finally, it is straightforward that $C[S_n,\theta]=\mb{S}^{n-1}$ for $\arccos\tfrac1n \le\theta\le\tfrac\pi2$. 
    \end{proof}
    
\begin{lemma}
    \label{lem:prob_comp}
Suppose $\psi^{\diamond}_{n}\in(\frac{\pi}{4}, \arccos\tfrac1{\sqrt{n}}]$ is arbitrary, and
\begin{equation*}
    B^{\diamond}_{n}(\theta):=2nA_n(\theta) - 2n(n-1)A_n\!\left(\theta,\tfrac\pi4\right), \quad\tfrac\pi4< \theta\le\tfrac\pi2. 
\end{equation*}
We have:
\begin{equation*}
\sigma(C[C_n,\theta]) \ 
\begin{cases}
=2nA_n(\theta), & 0<\theta\le \tfrac\pi4, \\
\ge B^{\diamond}_{n}(\theta), & \tfrac\pi4< \theta\le \psi^{\diamond}_{n}, \\
\ge B^{\diamond}_{n}(\psi^{\diamond}_{n}), & \psi^{\diamond}_{n}<\theta\le\arccos\!\tfrac{1}{\sqrt{n}}, \\
=1, & \arccos\!\tfrac{1}{\sqrt{n}}\leq \theta\le\tfrac\pi2.
\end{cases}
\end{equation*}
    \end{lemma}

    We omit the proof as it is completely similar to that of \cref{lem:Sprob_comp}.

\begin{remark}\label{rem:selection}
For each $n$, the parameters $\psi^{\scriptscriptstyle \triangle}_{n}$ and $\psi^{\diamond}_{n}$ are selected to numerically maximize $B^{\scriptscriptstyle \triangle}_{n}(\psi^{\scriptscriptstyle \triangle}_{n})$ and $B^{\diamond}_{n}(\psi^{\diamond}_{n})$, respectively.
\end{remark}

\subsection{A sphere packing bound}
The following lemma establishes that $n+2$ caps of radius $>\tfrac\pi4$ cannot be packed on $\mb{S}^{n-1}$, as shown in~\cite{Ra}*{Th.~1~(iii)}. Alternatively, a direct proof for $n=3$ that extends to higher dimensions can be found in the proof of Lemma~5 in~\cite{AKP-3d}.

\begin{lemma} 
\label{lemma:caps_pi/4} For any $\{x_1,\dots,x_{n+2}\}\subset\mb{S}^{n-1}$ there exist $1\le i<j\le n+2$ with $\theta(x_i,x_j)\le\tfrac\pi2$.
\end{lemma}


\section{Proof of \cref{thm:small} (dimensions $4\le n\le 15$)} \label{sec:small}

Suppose $K$ is a cap body in $\mb{R}^n$ with vertices $x_1, x_2, \dots, x_m$, and the corresponding radii of the caps $\varphi_1, \varphi_2,\dots, \varphi_m$. Without loss of generality, assume $\varphi_1\geq \varphi_2 \geq \ldots \geq \varphi_m$. 

To illuminate $K$, we start with a set of directions $\mc{T}$ which is the union of $s$ random rotations of $S_n$ (vertices of regular simplex) and $l$ random rotations of $C_n$ (vertices of cross-polytope), where $s$ and $l$, satisfying $s+l\ge 1$, will be selected later. It is routine to show that any rotation of either $S_n$ or $C_n$ satisfies condition~(ii) of~\cref{prop:ABPR}, and thus we will only need to verify~(i) of~\cref{prop:ABPR}. 

Since, by~\cref{lem:Sprob_comp} and~\cref{lem:prob_comp}, $C[S_n,\theta]=\mb{S}^{n-1}$ when $\theta \ge \arccos\left({\frac{1}{n}}\right)$, and $C[C_n,\theta]=\mb{S}^{n-1}$ when $\theta \ge \arccos\left({\frac{1}{\sqrt{n}}}\right)$, by \cref{prop:ABPR}~(i), any vertex $x_j$ of $K$ with $\varphi_j<\frac\pi2-\arccos\left({\frac{1}{n}}\right)$ is illuminated by one of the directions from $\mc{T}$. Moreover, if $\varphi_j=\frac\pi2-\arccos\left({\frac{1}{n}}\right)$, the probability that vertex $x_j$ is not illuminated is equal to 0 (vertex $x_j$ is not illuminated only when $l=0$ and for all $s$ random simplices $S_n$, the centre of one of the facets of $S_n$ has the same direction as $\hat{x}_j$). Thus, in what follows, we assume that $\varphi_m >  \frac\pi2-\arccos\frac{1}{n}$.

We will consider cases depending on the number of caps of similar size. Set $a_0:= \frac\pi2-\arccos\left({\frac{1}{n}}\right) \le\phi_m$. For a suitable positive integer $t$ that will be selected later, we define a discretization array $a=[a_0, a_1, \dots, a_t=\pi/2]$ where $a_i=a_0+\frac{a_t-a_0}{t}i$ for $0\leq i\leq t$. For any $1\leq j\leq m$ we have $\phi_j\in(a_0,a_t)$. Let $n_i$, $0\le i\le \mmone$, denote the number of indices $j$ such that $\varphi_j \in(a_i, a_{i+1}]$. The non-negative integers $n_i$, $0\le i\le \mmone$, will be the variables in the integer linear programming problem we are going to set up.

Next, let $E_{s,l}$ denote the expected number of caps not illuminated by $\mc{T}$. Since each such cap can be illuminated by an individual direction, we get the following bound on the illumination number: 
\begin{equation}\label{eqn:ill-bd0}
I(\mc{K}^n_c)\le \lfloor E_{s,l} \rfloor+(n+1)s+2nl.
\end{equation}
To conclude the proof, we will show that $E_{s,l}\le M_{s,l,t}$, where $M_{s,l,t}$ is the solution of a certain linear programming problem in $t$ variables.

More formally, note that by \cref{prop:ABPR}~(i) a vertex $x_j$ with $\varphi_j \in(a_i, a_{i+1}]$ is illuminated, provided at least one of the points of $\mc{T}$ is (strictly) within the geodesic distance $\tfrac\pi2-\phi_j$ of $\widehat{x}_j$. Therefore, the probability that $x_j$ is not illuminated by $\mc{T}$ is exactly
\[
\left(1-\sigma\left(C[S_n,{\tfrac\pi2-\phi_j}]\right)\right)^s\cdot \left(1-\sigma\left(C[C_n,{\tfrac\pi2-\phi_j}]\right)\right)^l,
\]
and so
\begin{equation}\label{eqn:target-prel}
E_{s,l}\le \sum_{i=0}^{t-1} n_i \left(1-\sigma\left(C[S_n,{\tfrac\pi2-a_{i+1}}]\right)\right)^s\cdot\left(1-\sigma\left(C[C_n,{\tfrac\pi2-a_{i+1}}]\right)\right)^l=:\sum_{i=0}^{t-1} n_i u_i^s v_i^l.
\end{equation}
Evaluating precisely $u_i$ and $v_i$ may be difficult (e.g. for certain values of angles three or more caps may overlap), but for our purposes, it suffices to upper bound $u_i$ and $v_i$. 
For this, we use \cref{lem:Sprob_comp,lem:prob_comp} with the parameters $\psi^{\scriptscriptstyle \triangle}_{n}$ and $\psi^{\diamond}_{n}$ chosen to numerically maximize $B^{\scriptscriptstyle \triangle}_{n}(\psi^{\scriptscriptstyle \triangle}_{n})$ and $B^{\diamond}_{n}(\psi^{\diamond}_{n})$, respectively. 
The estimates of $A_n(\theta)$ and $A_n(\alpha,\beta)$ are obtained on a computer using~\eqref{eqn:area-inter-equal}. (Implementation details will be discussed at the end of the proof.) Note that the lemmas provide lower bounds on $\sigma(C[S_n,\theta])$ or $\sigma(C[C_n,\theta])$ yielding upper bounds on $u_i$ and $v_i$. Denoting the corresponding calculated upper bounds $\ol{u}_i$ and $\ol{v}_i$, respectively, we obtain that $E_{s,l}$ 
does not exceed
\begin{equation}
    \label{eqn:target}
    \sum_{i=0}^{t-1} n_i (\ol u_i)^s (\ol v_i)^l,
\end{equation}
which will be the objective function in our integer programming problem. 

Besides non-negativity of $n_i$, we will use two constraints. First, recalling~\eqref{eqn:caps-not-overlap}, we get that the total measure of these caps is at most $1$. In terms of $n_i$, this provides the following constraint:
\begin{equation}
    \label{eqn:const bases pack}
    \sum_{i=0}^{t-1} n_iA_n(a_i)\le 1.
\end{equation}
The second constraint directly follows from~\eqref{eqn:caps-not-overlap} and \cref{lemma:caps_pi/4}:
\begin{equation}
    \label{eqn:const large caps}
    \sum_{0\le i<\m\,:\, a_i\ge \tfrac\pi4}n_i\le n+1.
\end{equation}

Let $M_{s,l,t}$ denote the solution of the integer linear programming problem 
    \begin{equation*}
        \text{\it maximize~\eqref{eqn:target} subject to~\eqref{eqn:const bases pack} and~\eqref{eqn:const large caps}}
    \end{equation*}
    with non-negative integer variables $n_i$, $0\le i\le t-1$, where in~\eqref{eqn:const bases pack} the coefficient $A_n(a_i)$ is evaluated numerically. By~\eqref{eqn:target-prel} we have $E_{s,l}\le M_{s,l,t}$ and so~\eqref{eqn:ill-bd0} implies 
    \begin{equation}\label{eqn:ill-bd}
    I(\mc{K}^n_c)\le \lfloor M_{s,l,t} \rfloor+(n+1)s+2nl.
    \end{equation}

Next, we establish an upper bound for $M_{s,l,t}$ via a computer-assisted 
continuous relaxation. We construct a modified linear program that strictly 
upper bounds $M_{s,l,t}$ by replacing the objective coefficients 
in~\eqref{eqn:target} with rational upper bounds, and the constraint 
coefficients in~\eqref{eqn:const bases pack} and \eqref{eqn:const large caps} 
with positive rational lower bounds. By relaxing the integer constraint on $n_i$, 
the resulting continuous linear program can be solved without floating-point 
roundoff errors using the exact rational arithmetic of the Parma Polyhedra 
Library (PPL) solver.

To obtain these rigorous rational bounds for the coefficients, we require 
mathematically certified estimates of $A_n(\theta)$ and $A_n(\alpha,\beta)$, 
which depend on the integral in~\eqref{eqn:area-inter-equal}. We evaluated 
these utilizing the Arb library, which employs arbitrary-precision complex ball 
arithmetic to provide guaranteed enclosures for all operations. Because Arb's 
rigorous integration algorithm requires the integrand to be analytic within a 
complex neighborhood of the integration interval, we applied the substitution 
$r=\sin(\theta)$ to eliminate fractional powers. The working precision for 
these computations was set to 64 bits.

We note that any selection of $\psi^{\scriptscriptstyle \triangle}_{n}$ and 
$\psi^{\diamond}_{n}$ yields a valid lower bound in 
\cref{lem:Sprob_comp,lem:prob_comp}. Consequently, computing strictly optimal 
values is unnecessary, and we utilized standard numerical precision to determine 
sufficiently tight approximations for $\psi^{\scriptscriptstyle \triangle}_{n}$ 
and $\psi^{\diamond}_{n}$.

All computations were executed via the SageMath (\cite{sagemath}) script 
provided in the appendix, requiring under an hour on a standard modern personal 
computer. For each dimension $n$ ($4 \le n \le 15$), we selected specific 
parameters $t$ ($50 \le t \le 900$), $s$, and $l$ to minimize the resulting 
bound on $I(\cK_c^n)$ via~\eqref{eqn:ill-bd}. A summary of these results is 
presented in \cref{tbl}.

\begin{table}[h!]
	\begin{center}
		\begin{tabular}{|c|c|c|c|}
			\hline
			$n$ & $I(\mc{K}^n_c)\le$ & $s$ (simplices) & $l$ (cross-polytopes)  \\
			\hline
            \hline
			$4$ & $11$ & $1$ & $0$ \\
			\hline
			$5$ & $17$ & $0$ & $1$ \\
			\hline
			$6$ & $27$ & $0$ & $1$ \\
			\hline
			$7$ & $44$ & $1$ & $1$ \\
			\hline
			$8$ & $66$ & $0$ & $2$ \\
			\hline
			$9$ & $103$ & $0$ & $3$ \\
			\hline
			$10$ & $157$ & $0$ & $4$ \\
			\hline
			$11$ & $236$ & $0$ & $6$ \\
			\hline
			$12$ & $353$ & $0$ & $8$ \\
			\hline
			$13$ & $520$ & $0$ & $11$ \\
            \hline
			$14$ & $769$ & $0$ & $15$ \\
            \hline
			$15$ & $1127$ & $0$ & $21$ \\
			\hline
		\end{tabular}
	\end{center}
	\caption{Upper bounds on $I(\mathcal{K}_{c}^{n})$ for $4\le n \le 15$.}\label{tbl}
\end{table}


\section{Proof of \cref{thm:explicit} (dimensions $n\ge 9$)} \label{sec:large}

Assume that $K$ is a cap body in $\mm^n$ for $n\geq 9$. Once again (as in Section~\ref{sec:small}), assume $x_1, \ldots, x_m$ are vertices of $K$ and the corresponding radii of the caps are $\phi_1, \ldots, \phi_{m}$. 

We will begin by illuminating $K$ with the union of $x\ge 1$ independent random rotations of the cross-polytope $C_n$. It is immediate that any rotation of $C_n$ satisfies \cref{prop:ABPR}~(ii), so it remains only
to verify \cref{prop:ABPR}~(i). Also, by Lemma~\ref{lem:prob_comp}, $
\sigma(C[C_n,\theta])=1$ when $\theta \ge \arccos{\frac{1}{\sqrt{n}}}$. Therefore, any cap with radius smaller than $\frac{\pi}{2}-\arccos{\frac{1}{\sqrt{n}}}$ will be illuminated by any rotation of $C_n$. On the other hand, by \cref{lemma:caps_pi/4}, there are at most $n+1$ caps of radius $>\pi/4$.

Now let $m^\prime$ be the number of indices $j$ such that $\phi_j\in\left[\frac{\pi}{2}-\arccos{\frac{1}{\sqrt{n}}}, \frac{\pi}{4} \right]$. Define $p:=A_n\left(\frac{\pi}{4}\right)$, and $q:=A_n\left(\frac{\pi}{2}-\arccos{\frac{1}{\sqrt{n}}}\right)$. Note that the $2n$ caps in $C[C_n,\pi/4]$ have pairwise disjoint relative interiors, so $2np<1$. A vertex $x_j$ of $K$ is not illuminated by any of $x$ rotations of $C_n$ with probability at most $(1-2np)^x$, so we can choose the rotations in a way that there are at most $(1-2np)^xm^\prime$ unilluminated vertices of $K$. All together, taking into account that $m^\prime q\le 1$ as base caps do not overlap, we see that $K$ can be illuminated by at most
\begin{equation}
    (1-2np)^xm^\prime + 2nx + n+1\leq (1-2np)^x\tfrac1q + 2nx + n+1=:f(x)  
\end{equation}
directions. Now we want to minimize $f(x)$ over $x\in \mb{Z}^+$.

Note that $f'(x)=(1-2np)^x\frac{1}{q}\ln{(1-2np)} + 2n$, and the only point $x_0>0$ satisfying $f'(x_0)=0$ is
\begin{equation}
    \label{eqn:x0solution}
    x_0=\frac1{-\ln(1-2np)}\left(\ln(-\ln(1-2np))-\ln(2nq)\right).
\end{equation}
Thus, $\min_{x\in\mb{Z}^+}f(x)=f(y_0)$ must be attained at a point $y_0=\lfloor x_0\rfloor$ or $y_0=\lceil x_0 \rceil$. We get \begin{equation}\label{eqn:MVT}
    f(y_0)= f'(z_0)(y_0-x_0)+f(x_0),
\end{equation}
where $|z_0-x_0|< 1$. Now, suppose that $z_0=x_0+h$, then by $f'(x_0)=0$
\begin{align*}
    f'(z_0)&= f'(x_0+h)\\
    &= (1-2np)^h(1-2np)^{x_0}\frac{1}{q}\ln{(1-2np)} + 2n\\
    &= 2n(1-(1-2np)^h)=:g(h).
\end{align*}
Since $g$ is a monotone function of $h\in[-1,1]$, we have $\max_{h\in[-1,1]}|g(h)|$ is achieved either at $h=1$ or $h=-1$, and a direct verification gives $|g(-1)|>|g(1)|$.
%
%
We conclude that 
$|f'(z_0)|\le |g(-1)|=\frac{4n^2p}{1-2np}$. By~\eqref{eqn:MVT}, $|f(y_0)-f(x_0)|\le |f'(z_0)|$, so with the help of~\eqref{eqn:x0solution}, we get 
\begin{align}
    f(y_0)&\le \frac{4n^2p}{1-2np} + f(x_0) \nonumber\\
    &= \frac{4n^2p}{1-2np} + \frac{2n}{-\ln{(1-2np)}} + \frac{2n}{-\ln(1-2np)}\left(\ln(-\ln(1-2np))-\ln(2nq)\right) + n+1. \label{eqn:yi}
\end{align}

Using (\ref{eqn:BW}) with angles $\theta=\pi/4, \pi/2-\arccos(1/\sqrt{n})$ we obtain estimates on $p,q$
\begin{equation}\label{eqn:pq}
    \left(\frac{1}{\sqrt{2}}\right)^n\frac{1}{\sqrt{\pi n}}\leq p \leq \left(\frac{1}{\sqrt{2}}\right)^n\sqrt{\frac{2}{\pi (n-1)}}, \qquad    \frac{1}{\sqrt{2 \pi}}\left(\frac{1}{\sqrt{n}}\right)^n \le q.
\end{equation}
We will use the estimates $x/(1+x)\le \ln(1+x)\le x $ which are valid for $x>-1$. We have $-\ln(1-2np)\le 2np/(1-2np)\le 3/n$ for all $n\geq 9$, where the second inequality is verified by direct computation for initial values and standard calculus for larger values. Finally, the last inequality also implies that $\frac{4n^2p}{1-2np}\leq 6$ for all $n\geq 9$.

So, continuing~\eqref{eqn:yi},
\begin{align*}
    f(y_0)  
    &\le 6+\frac{2n}{2np}+\frac{2n}{2np}\left(\ln\left(-\frac{-2np}{1-2np}\right)-\ln(2nq))\right)+n+1\\
    &= n+7+\frac{1}{p}\left(1+\ln(p/q)-\ln(1-2np)\right)\\
    &\leq n+7 +\frac{1}{p}(1+\ln(p/q)+3/n)\\
    & \leq n+7+(\sqrt{2})^n\sqrt{\pi n}\left(1+3/n+\ln\left((\sqrt{n/2})^n \frac{2}{\sqrt{n-1}}\right)\right)\\
    &\leq n+7+(\sqrt{2})^n\sqrt{\pi n}\left(1+\frac3n+\frac12n\ln\frac n2 +\ln\frac2{\sqrt{n-1}}\right),\\
\end{align*}
which is the claimed bound for $n\ge9$. It is routine to check that this bound does not exceed $2^n-1$ for $n\geq 13$.

\section{Discussion and remarks}\label{sec:remarks}

For each dimension $n$ and fixed $s$, $l$, the best bound one can obtain with our method is $I(\mc{K}^n_c)\le \liminf_{t\to\infty}\lfloor M_{s,l,t} \rfloor+(n+1)s+2nl$. 
We varied the number of intervals $t$ between $50$ and $900$ selecting a smaller value when possible which still yields the same upper bound on $I(\mc{K}_c^n)$. Values of $s$ and $l$ were chosen on a case-by-case basis to fit the dimension $n$. It is likely that the resulting bounds we obtained for $4\le n\le 8$ are best possible by the method, in other words, taking larger $t$ or other $s$ and $l$ will not yield any improvement. For other dimensions, one can further increase $t$ at the expense of longer computations and improve the estimates in \cref{tbl}. We note, however, that the computation cost of the linear program may be non-linear in $t$ and could increase rapidly.

If one is interested only in establishing $I(\mc{K}_c^n)<2^n$ for $4\le n\le 15$, then $t$ can be chosen to be much smaller, say $t=8$. The corresponding values of $s$, $l$ are available in the last portion of the script.

It is interesting to note that for $n=4$ using a (random rotation of) simplex is better, while for $n=5,6$ using a cross-polytope is better. Confirming $I(\mc{K}_c^n)<2^n$ could have been done using only simplices or using only cross-polytopes, but we chose to obtain better bounds and include both. It is apparent from \cref{tbl}, where we found heuristically optimal choices of $s$ and $l$, that while using only cross-polytopes is often optimal among the choices considered, both configurations are useful for obtaining better estimates for $n=4,7$. One can also try random rotations of other ``good" point configurations available in specific dimensions.

Finally, the definition of cap bodies can be modified to allow infinitely many vertices, following the approach in~\cite{IS}. Under this modification, \cref{prop:ABPR} remains valid, and there are only finitely many caps with a radius exceeding any given positive value. All caps with sufficiently small radii are illuminated by an arbitrary rotation of a simplex or cross-polytope, while the treatment of larger caps presented in this work still applies. Consequently, our results extend to these infinite-vertex cap bodies.

\section{Acknowledgments}
We thank the referee for a very careful reading of our work and numerous comments, and, in particular, for identifying an oversight in the first version of this manuscript where the selections of $\psi^{\scriptscriptstyle \triangle}_{n}$ and $\psi^{\diamond}_{n}$ for applications of \cref{lem:Sprob_comp,lem:prob_comp} were not optimal. 
Addressing this ultimately led to improved bounds on the illumination number in dimensions $n=6$ and $n\ge 8$.

\section*{Appendix}

\noindent Computation of $\psi^{\scriptscriptstyle \triangle}_{n}$ and $\psi^{\diamond}_{n}$ as per \cref{rem:selection}:
\begin{python}
#numerical integration is used for the purpose of pre-computing optimal psi values
def ternary_search_max(func, n_val, left, right, max_iter):
    for i in range(max_iter):
        m1 = left + (right - left) / 3.0
        m2 = right - (right - left) / 3.0
        f_m1 = func(n_val, m1)
        f_m2 = func(n_val, m2)
        if f_m1 < f_m2:
            left = m1
        else:
            right = m2
    return (left + right) / 2.0

def inter(n,alpha,beta,mult=0): #probabilistic measure of the intersection of two caps with centres 2*beta apart, radius alpha as per (4), plus the error of numerical integration multiplied by mult
    integ = numerical_integral(lambda r: (n-2)/pi*(1-r^2)^((n-4)/2)*(arccos(cos(alpha)/r)-beta)*r, cos(alpha)/cos(beta), 1)
    return integ[0]+mult*integ[1]

def cross_cap_prob(n,theta): #probabilistic measure of the union of spherical caps centred at cross polytopes, for theta>=pi/4, 
    #when theta is too large, this will return a lower bound
    return 2*n*inter(n,theta,0,-1)-2*n*(n-1)*inter(n,theta,pi/4,1)

def simplex(n,theta): #probabilistic measure of the union of spherical caps centred at simplex
    return (n+1)*inter(n,theta,0,-1)-1/2*n*(n+1)*inter(n,theta,arccos(-1/n)/2,1)

maxt_c = [0]*16
maxt_s = [0]*16
for n in range(4,16):
    maxt_c[n] = ternary_search_max(cross_cap_prob, n, (pi/4).n(), (arccos(1/sqrt(n))).n(), 20)
    maxt_s[n] = ternary_search_max(simplex, n, (1/2 * arccos(-1/n)).n(), (arccos(1/n)).n(), 20)

print("maxt_c = "+str(maxt_c),"maxt_s = "+str(maxt_s))
\end{python}

\vspace{2.5 mm}

\noindent Main script:
\begin{python}
from sage.all import *
#pre-computed values of psi
maxt_c = [0, 0, 0, 0, 0.982801897660807, 1.00524007590983, 1.02421744871292, 1.04064922651209, 1.05510104120268, 1.06796254888841, 1.07961984845650, 1.09013610877061, 1.09984487751988, 1.10879399390152, 1.11711030431677, 1.12481703731639]
maxt_s = [0, 0, 0, 0, 1.20127138942692, 1.18892069741706, 1.18418710112405, 1.18325482143378, 1.18444710830655, 1.18678299779178, 1.18984586047734, 1.19321483347841, 1.19684205678909, 1.20051689748862, 1.20426109818127, 1.20794267994429]

def inter_certified(n, alpha_val, beta_val, mult=-1, prec=64):
    """
    Computes a mathematically certified bound for the integral in (4) using Arb (ComplexBallField).
    mult = -1 returns the exact rational lower bound.
    mult = 1 returns the exact rational upper bound.
    """
    CBF = ComplexBallField(prec)
    alpha = CBF(alpha_val)
    beta = CBF(beta_val)
    
    # Calculate limits
    lower_lim_ratio = alpha.cos() / beta.cos()
    
    # Strictly check if the lower limit exceeds 1 
    if lower_lim_ratio.real().endpoints()[0] >= 1:
        return QQ(0)
    
    # To prevent domain errors from minor interval overshoot, cap the ratio at 1
    if lower_lim_ratio.real().endpoints()[1] > 1:
        lower_lim_ratio = CBF(1)
        
    theta_min = lower_lim_ratio.arcsin()
    theta_max = CBF(pi)/2
    
    # The analytically smoothed integrand using r = sin(theta)
    def integrand(theta, _):
        term1 = CBF(n-2) / CBF(pi)
        term2 = theta.cos()**(n-3) * theta.sin()
        term3 = (alpha.cos() / theta.sin()).arccos() - beta
        return term1 * term2 * term3
        
    # Rigorous integration
    res = CBF.integral(integrand, theta_min, theta_max)
    
    # Extract strict endpoints as exact rationals
    lower_bound_exact = QQ(res.real().endpoints()[0])
    upper_bound_exact = QQ(res.real().endpoints()[1])
    
    if mult == -1:
        return lower_bound_exact
    elif mult == 1:
        return upper_bound_exact

def cross_cap_prob_certified(n, theta): #rational lower bound on sigma(C[C_n,theta])
    term1 = 2 * n * inter_certified(n, theta, 0, -1)
    term2 = 2 * n * (n - 1) * inter_certified(n, theta, pi/4, 1)
    return term1 - term2

def simplex_certified(n, theta): #rational lower bound on sigma(C[S_n,theta])
    term1 = (n + 1) * inter_certified(n, theta, 0, -1)
    term2 = QQ(n * (n + 1)/2) * inter_certified(n, theta, arccos(-1/n)/2, 1)
    return term1 - term2

def min_illum_certified(n, t, s, c):
    A = pi/2 - arccos(1/n)
    B = pi/2
    a = [A + (B - A)*i/t for i in range(t+1)] 
    
    probc = [] #rational lower bounds on sigma(C[C_n,pi/2-a[i]])
    probs = [] #rational lower bounds on sigma(C[S_n,pi/2-a[i]])
    
    pc_max = cross_cap_prob_certified(n, maxt_c[n]) #lower bound
    
    for i in range(1, t+1): #note loop starts with 1, so probc and probs indexes are smaller by 1
        alpha = a[i]
        if alpha <= pi/2 - arccos(1/sqrt(n)): #covering radius for cross polytope
            probc.append(QQ(1))
        elif alpha <= pi/2 - maxt_c[n]: #to the right of max to optimize use of Bonferroni
            probc.append(pc_max)
        elif alpha < pi/4: #double/Bonferroni
            probc.append(cross_cap_prob_certified(n, pi/2 - alpha))
        else:  #no intersections
            probc.append(2 * n * inter_certified(n, pi/2 - alpha, 0, -1))
            
    ps_max = simplex_certified(n, maxt_s[n]) #lower bound
    
    for i in range(1, t+1):
        alpha = a[i]
        if alpha < pi/2 - arccos(1/n): #never true, for homogeneity
            probs.append(QQ(1))
        elif alpha <= pi/2 - maxt_s[n]: #to the right of max to optimize use of Bonferroni
            probs.append(ps_max)
        elif alpha < pi/2 - 1/2*arccos(-1/n): #double/Bonferroni
            probs.append(min(QQ(1), simplex_certified(n, pi/2 - alpha)))
        else:  #no intersections
            probs.append((n + 1) * inter_certified(n, pi/2 - alpha, 0, -1))
            
    #exact rational linear programming (variables allowed to be non-integer) using the 'ppl' solver
    p = MixedIntegerLinearProgram(solver='ppl')
    v = p.new_variable(nonnegative=True)
    
    obj_expr = 0
    for i in range(t):
        obj_coeff = (QQ(1) - probc[i])**c * (QQ(1) - probs[i])**s #this is an upper bound on the expected number of non-illuminated caps, because probc[i] and probs[i] are lower bounds on the probability
        obj_expr += v[i] * obj_coeff
        
    p.set_objective(obj_expr) 
    
    con1_expr = 0
    for i in range(t):
        con_coeff = inter_certified(n, a[i], 0, -1) #constraint coefficients use lower bounds to enlarge the feasible set
        con1_expr += v[i] * con_coeff
        
    p.add_constraint(con1_expr <= QQ(1))
    
    con2_expr = 0
    for i in range(t):
        if a[i] >= pi/4:
            con2_expr += v[i] * QQ(1)
            
    p.add_constraint(con2_expr <= QQ(n + 1))
    
    exact_max = p.solve()
    
    # converting the exact result to a decimal representation for readable terminal output
    target = round(float(exact_max) + 2*n*c + s*(n+1), 4)
    return n, target

print(min_illum_certified(4, 50, 1, 0))
print(min_illum_certified(5, 100, 0, 1))
print(min_illum_certified(6, 100, 0, 1))
print(min_illum_certified(7, 200, 1, 1))
print(min_illum_certified(8, 700, 0, 2))
print(min_illum_certified(9, 400, 0, 3))
print(min_illum_certified(10, 400, 0, 4))
print(min_illum_certified(11, 700, 0, 6))
print(min_illum_certified(12, 500, 0, 8))
print(min_illum_certified(13, 850, 0, 11)) 
print(min_illum_certified(14, 700, 0, 15)) 
print(min_illum_certified(15, 900, 0, 21)) 

#smaller t (=8) only to confirm conjecture
print(min_illum_certified(4,8,0,1),2^4)
print(min_illum_certified(5,8,0,1),2^5)
print(min_illum_certified(6,8,0,2),2^6)
print(min_illum_certified(7,8,0,3),2^7)
print(min_illum_certified(8,8,0,3),2^8)
print(min_illum_certified(9,8,0,6),2^9)
print(min_illum_certified(10,8,0,11),2^10)
print(min_illum_certified(11,8,0,18),2^11)
print(min_illum_certified(12,8,0,29),2^12)
print(min_illum_certified(13,8,0,45),2^13)
print(min_illum_certified(14,8,0,68),2^14)
print(min_illum_certified(15,8,0,100),2^15)
\end{python}

\vspace{2.5 mm}

\noindent Output:
\begin{python}
(4, 11.7639)
(5, 17.8327)
(6, 27.9351)
(7, 44.6747)
(8, 66.9327)
(9, 103.9368)
(10, 157.9018)
(11, 236.592)
(12, 353.7566)
(13, 520.6181)
(14, 769.2394)
(15, 1127.3994)

(4, 15.0119) 16
(5, 25.6834) 32
(6, 48.3176) 64
(7, 85.0816) 128
(8, 145.9289) 256
(9, 251.8416) 512
(10, 426.103) 1024
(11, 708.0317) 2048
(12, 1158.9933) 4096
(13, 1872.7889) 8192
(14, 2995.3914) 16384
(15, 4759.3239) 32768
\end{python}

\end{document}